\documentclass[10pt]{article}
\font\smallit=cmti10
\font\smalltt=cmtt10

\usepackage{amssymb,latexsym,amsmath,epsfig,amsthm,url} 

\makeatletter

\renewcommand\section{\@startsection {section}{1}{\z@}
{-30pt \@plus -1ex \@minus -.2ex}
{2.3ex \@plus.2ex}
{\normalfont\normalsize\bfseries}}

\renewcommand\subsection{\@startsection{subsection}{2}{\z@}
{-3.25ex\@plus -1ex \@minus -.2ex}
{1.5ex \@plus .2ex}
{\normalfont\normalsize\bfseries}}

\renewcommand{\@seccntformat}[1]{\csname the#1\endcsname. }

\makeatother

\newtheorem{theorem}{Theorem}
\newtheorem{lemma}{Lemma}

\newtheorem{corollary}{Corollary}

\newtheoremstyle{defi}
  {10pt}		  
  {10pt}  
  {\rm}  
  {\parindent}     
  {\bf}  
  {. }    
  { }    
  {}     
\theoremstyle{defi}

\newtheorem{remark}[theorem]{Remark}

\begin{document}

\begin{center}
\uppercase{\bf The Non-Euler Part of a Spoof Odd Perfect Number is Not Almost Perfect}
\vskip 20pt
{\bf Jose Arnaldo B. Dris\footnote{PhD Student, University of the Philippines - Diliman}}\\
{\smallit Institute of Mathematics, University of the Philippines \\ 
Diliman, Quezon City}\\
{\tt josearnaldobdris@gmail.com}, {\tt jadris@feu.edu.ph}\\ 
\end{center}
\vskip 30pt

\centerline{\bf Abstract}
\noindent
We call $n$ a spoof odd perfect number if $n$ is odd and $n=km$ for two integers $k,m>1$ such that $\sigma(k)(m+1)=2n$, where $\sigma$ is the sum-of-divisors function.  In this paper, we show how results analogous to those of odd perfect numbers could be established for spoof odd perfect numbers (otherwise known in the literature as Descartes numbers).  In particular, we prove that $k$ is not almost perfect.

\pagestyle{myheadings} 
\markright{\smalltt \hfill} 
\thispagestyle{empty} 
\baselineskip=12.875pt 
\vskip 30pt

\section{Introduction}
This article is an elucidation of some of the recent discoveries/advances in the preprint \cite{1}, as applied to the case of spoof odd perfect numbers, with details in the older version of this preprint titled \emph{The Abundancy Index of Divisors of Spoof Odd Perfect Numbers}.

Recall that we call $n$ a spoof odd perfect number if $n$ is odd and $n=km$ for two integers $k, m > 1$, such that $\sigma(k)(m+1)=2n=2km$.  The number $m$ is called the \emph{quasi-Euler prime} of the spoof $n$, while $k$ is referred to as the \emph{non-Euler part} of $n$.

For some recent papers on spoofs (otherwise known as Descartes numbers), we refer the interested reader to \cite{6} and \cite{7}.

\section{Preliminaries}
We begin with the following very useful lemmas.  (In what follows, we take $\sigma(x)$ to be the sum of the divisors of $x$, and denote the abundancy index of $x \in \mathbb{N}$ as $I(x)=\sigma(x)/x$, where $\mathbb{N}$ is the set of natural numbers or positive integers.)

\begin{lemma}\label{lem:1} 
Let $n = km$ be a spoof odd perfect number.  Then
$$\frac{\sigma(m)}{\sqrt{k}}+\frac{\sigma(\sqrt{k})}{m} \text{ is bounded } \iff \frac{m}{\sqrt{k}}+\frac{\sqrt{k}}{m} \text{ is bounded }.$$
\end{lemma}

\begin{remark}\label{rem:1}
Note that, in Lemma \ref{lem:1}, we define
$$\sigma(m) := m+1.$$
\end{remark}

\begin{proof}
The claimed result follows from the inequalities
$$m < \sigma(m) = m + 1 < 2m$$
$$\sqrt{k} < \sigma(\sqrt{k}) < 2\sqrt{k}$$
(since $\sqrt{k}$ is a proper divisor of the deficient number $k$). Consequently,
$$\frac{m}{\sqrt{k}} < \frac{\sigma(m)}{\sqrt{k}} < 2\cdot\frac{m}{\sqrt{k}}$$
$$\frac{\sqrt{k}}{m} < \frac{\sigma(\sqrt{k})}{m} < 2\cdot\frac{\sqrt{k}}{m}$$
from which it follows that
$$\frac{m}{\sqrt{k}}+\frac{\sqrt{k}}{m} < \frac{\sigma(m)}{\sqrt{k}}+\frac{\sigma(\sqrt{k})}{m} < 2\cdot\bigg(\frac{m}{\sqrt{k}}+\frac{\sqrt{k}}{m}\bigg).$$
We therefore conclude that
$$\frac{\sigma(m)}{\sqrt{k}}+\frac{\sigma(\sqrt{k})}{m} \text{ is bounded } \iff \frac{m}{\sqrt{k}}+\frac{\sqrt{k}}{m} \text{ is bounded },$$
as desired.
\end{proof}

\begin{remark}\label{rem:2}
In general, the function
$$f(z) := z + \frac{1}{z}$$
is not bounded from above.  (It suffices to consider the cases $z \to 0^{+}$ and $z \to +\infty$.)  This means that we do not expect the sum
$$\frac{\sigma(m)}{\sqrt{k}}+\frac{\sigma(\sqrt{k})}{m}$$
to be bounded from above.
\end{remark}

\begin{corollary}\label{cor:1}
Let $n = km$ be a spoof odd perfect number.  Then the following conditions hold:
\begin{enumerate}
{
\item{$\sigma(\sqrt{k}) \neq m+1=\sigma(m)$}
\item{$\sigma(\sqrt{k}) \neq m$}
}
\end{enumerate}
\end{corollary}

\begin{proof}
(a)  Suppose that $n=km$ is a spoof satisfying the condition
$$\sigma(\sqrt{k}) = \sigma(m).$$
It follows that
$$\frac{\sigma(\sqrt{k})}{\sqrt{k}} = \frac{\sigma(m)}{\sqrt{k}}$$
and
$$\frac{\sigma(\sqrt{k})}{m} = \frac{\sigma(m)}{m}$$
from which we obtain
$$\frac{\sigma(m)}{\sqrt{k}}+\frac{\sigma(\sqrt{k})}{m}=\frac{\sigma(\sqrt{k})}{\sqrt{k}}+\frac{\sigma(m)}{m}.$$
But we also have
$$\frac{\sigma(\sqrt{k})}{\sqrt{k}}+\frac{\sigma(m)}{m} < \frac{\sigma(k)}{k}+\frac{\sigma(m)}{m} = I(k) + \frac{m+1}{m} < 2 + \frac{10}{9} = \frac{28}{9}.$$
Finally, we get
$$\frac{\sigma(m)}{\sqrt{k}}+\frac{\sigma(\sqrt{k})}{m} < \frac{28}{9},$$
which contradicts Lemma \ref{lem:1} and Remark \ref{rem:2}.  We conclude that
$$\sigma(\sqrt{k}) \neq m+1=\sigma(m).$$

(b)  Suppose that $n=km$ is a spoof satisfying the condition
$$\sigma(\sqrt{k}) = m.$$
It follows that
$$\frac{\sigma(\sqrt{k})}{m}=1.$$
But
$$1 < \frac{\sigma(\sqrt{k})}{\sqrt{k}}\cdot\frac{\sigma(m)}{m} = \frac{\sigma(\sqrt{k})}{m}\cdot\frac{\sigma(m)}{\sqrt{k}} < 2.$$
This implies that
$$1 < \frac{\sigma(m)}{\sqrt{k}} < 2.$$
In particular,
$$\frac{\sigma(\sqrt{k})}{m}+\frac{\sigma(m)}{\sqrt{k}}=1+\frac{\sigma(m)}{\sqrt{k}} < 1+2 = 3.$$
This contradicts Lemma \ref{lem:1} and Remark \ref{rem:2}.  We conclude that
$$\sigma(\sqrt{k}) \neq m.$$
\end{proof}

\begin{lemma}\label{lem:2}
Let $a, b \in \mathbb{N}$.
\begin{enumerate}
{
\item{If $I(a) + I(b) < \sigma(a)/b + \sigma(b)/a$, then $a < b \iff \sigma(a) < \sigma(b)$ holds.}
\item{If $\sigma(a)/b + \sigma(b)/a < I(a) + I(b)$, then $a < b \iff \sigma(b) < \sigma(a)$ holds.}
\item{If $I(a) + I(b) = \sigma(a)/b + \sigma(b)/a$ holds, then either $a = b$ or $\sigma(a) = \sigma(b)$ is true.}
}
\end{enumerate}
\end{lemma}

\begin{proof}
We refer the interested reader to a proof of part (2) in page 6 of this preprint \cite{1}.  The proofs for parts (1) and (3) are very similar.
\end{proof}

\begin{remark}\label{rem:3}
Note that if we let $a = m$ and $b = \sqrt{k}$ in Lemma \ref{lem:2}, and if we make the additional assumption that $\gcd(a, b) = \gcd(m, k) = 1$, then case (3) is immediately ruled out, as $\gcd(m, \sqrt{k}) = 1$ implies that $m \neq \sqrt{k}$.  Additionally, note that $\sigma(m) \neq \sigma(\sqrt{k})$ per Corollary \ref{cor:1} (a).

Likewise, note that Lemma \ref{lem:1} and Remark \ref{rem:2} rules out case (2), as it implies that
$$\frac{\sigma(m)}{\sqrt{k}}+\frac{\sigma(\sqrt{k})}{m} < I(m) + I(\sqrt{k}) < \frac{m+1}{m} + I(k) < \frac{10}{9}+2=\frac{28}{9},$$
a contradiction.

Hence, we are left with the scenario under case (1):
$$I(m)+I(\sqrt{k}) = \frac{m+1}{m}+\frac{\sigma(\sqrt{k})}{\sqrt{k}} < \frac{m+1}{\sqrt{k}}+\frac{\sigma(\sqrt{k})}{m},$$
which per Lemma \ref{lem:2} implies that
$$m < \sqrt{k} \iff m+1 < \sigma(\sqrt{k}).$$
\end{remark}

The considerations in Remark \ref{rem:3} prove the following proposition.

\begin{theorem}\label{thm:1}
Let $n = km$ be a spoof odd perfect number.  Then the series of biconditionals
$$m < \sqrt{k} \iff m+1 < \sigma(\sqrt{k}) \iff \frac{m+1}{\sqrt{k}} < \frac{\sigma(\sqrt{k})}{m}$$
hold.
\end{theorem}

\begin{proof}
The proof is trivial.
\end{proof}

\begin{remark}\label{rem:4}
Note that
$$\frac{m+1}{\sqrt{k}} \neq \frac{\sigma(\sqrt{k})}{m}$$
is in general true under the assumption $\gcd(m, \sqrt{k})=1$, since it follows from the fact that
$$1 < \frac{m+1}{\sqrt{k}}\cdot\frac{\sigma(\sqrt{k})}{m}=\frac{m+1}{m}\cdot\frac{\sigma(\sqrt{k})}{\sqrt{k}} < 2.$$

Also, note that since $m \neq \sqrt{k}$ (because $\gcd(m, \sqrt{k})=1$), and $m+1 = \sigma(m) \neq \sigma(\sqrt{k})$ (by Corollary \ref{cor:1}, (a)), then equivalently, we have the series of biconditionals
$$\sqrt{k} < m \iff \sigma(\sqrt{k}) < m+1 \iff \frac{\sigma(\sqrt{k})}{m} < \frac{m+1}{\sqrt{k}}.$$
\end{remark}

\section{Main Results}

\begin{remark}\label{rem:5}
Let us double-check the findings of Theorem \ref{thm:1} (and Remark \ref{rem:4}) using the only spoof that we know of, as a test case.

In the Descartes spoof,
$$m = 22021 = {{19}^2}\cdot{61}$$
and
$$\sqrt{k} = {3}\cdot{7}\cdot{11}\cdot{13} = 3003$$
so that we obtain
$$m + 1 = 22022 = 2\cdot{11011}$$
and
$$\sigma(\sqrt{k}) = (3+1)\cdot(7+1)\cdot({11}+1)\cdot({13}+1) = 5376 = {2^8}\cdot{3}\cdot{7}.$$
Notice that we then have
$$\sqrt{k} < m,$$
$$\sigma(\sqrt{k}) < m+1,$$
and
$$\frac{\sigma(\sqrt{k})}{m} = \frac{5376}{22021} < 1 < \frac{22022}{3003} = \frac{m+1}{\sqrt{k}},$$
in perfect agreement with the results in Theorem 3.
\end{remark}

\begin{remark}\label{rem:6}
Using Theorem 3, we list down all allowable permutations of the set
$$\left\{m, m+1, \sqrt{k}, \sigma(\sqrt{k})\right\}$$
(following the usual ordering on $\mathbb{N}$) subject to the biconditional
$$m < \sqrt{k} \iff m+1 < \sigma(\sqrt{k}) \iff \frac{m+1}{\sqrt{k}} < \frac{\sigma(\sqrt{k})}{m}$$
and the constraints
$$m < m + 1, \text{   } \sqrt{k} < \sigma(\sqrt{k}).$$

We have the following cases to consider: \\

Case A:  $m < m+1 := \sigma(m) < \sqrt{k} < \sigma(\sqrt{k})$

Under this case, $m < \sqrt{k} < k$ which is true if and only if $k$ is not an odd almost perfect number (see page 6 of this preprint \cite{2}).  That is, 
$$D(k) := 2k - \sigma(k) = \frac{\sigma(k)}{m} > \frac{\sigma(\sqrt{k})}{m} > 1,$$
whence there is no contradiction.

Case B:  $\sqrt{k} < \sigma(\sqrt{k}) < m < m + 1$

Noting that, for the Descartes spoof, we have
$$D(k) = D({3003}^2) = 2(3003)^2 - \sigma({3003}^2) = 819 > 1,$$
so that $k$ is not almost perfect, it would seem prudent to try to establish the inequality $m < k$ for Case B.

To do so, one would need to repeat the methodology starting with Lemma \ref{lem:2}, but this time with $a = m$ and $b = k$.  Again, we are left with the scenario under case (1):
$$I(m)+I(k) = \frac{m+1}{m}+\frac{\sigma(k)}{k} < \frac{m+1}{k}+\frac{\sigma(k)}{m},$$
which per Lemma \ref{lem:2} implies that
$$m < k \iff m+1 < \sigma(k).$$
\end{remark}

Thus we see that we have the following corollary to Theorem \ref{thm:1}.

\begin{corollary}\label{cor:2}
Let $n = km$ be a spoof odd perfect number.  Then the series of biconditionals
$$m < k \iff m+1 < \sigma(k) \iff \frac{m+1}{k} < \frac{\sigma(k)}{m}$$
hold.
\end{corollary}

\begin{remark}\label{rem:7}
Note that
$$\frac{m+1}{k} \neq \frac{\sigma(k)}{m}$$
is in general true under the assumption $\gcd(m, k)=1$, since it follows from the fact that
$$1 < \frac{m+1}{k}\cdot\frac{\sigma(k)}{m}=\frac{m+1}{m}\cdot\frac{\sigma(k)}{k} = 2.$$
Otherwise, we have $2 = (m+1)/k$ (which implies that $m=2k-1$) and $\sigma(k)/m=1$, so that $D(k) = 2k - \sigma(k) = 1$, which is equivalent to $k < m$.  (This is not true under Case A above, so it suffices to prove directly that $m < k$ under Case B above.)

Also, note that since $m \neq k$ (because $\gcd(m, k)=1$), and $m+1 = \sigma(m) \neq \sigma(k)$ (by a similar result as proved in Corollary \ref{cor:1}, (a)), then equivalently, we have the series of biconditionals
$$k < m \iff \sigma(k) < m+1 \iff \frac{\sigma(k)}{m} < \frac{m+1}{k}.$$
Lastly, note that $\sigma(k) \neq m$ (by a similar result as proved in Corollary \ref{cor:1}, (b)).
\end{remark}

\begin{remark}\label{rem:8}
Using Corollary \ref{cor:2} (and the considerations in Remark \ref{rem:7}), we list down all allowable permutations of the set
$$\left\{m, m+1, k, \sigma(k)\right\}$$
(following the usual ordering on $\mathbb{N}$) subject to the biconditional
$$m < k \iff m+1 < \sigma(k) \iff \frac{m+1}{k} < \frac{\sigma(k)}{m}$$
and the constraints
$$m < m + 1, \text{   } k < \sigma(k).$$

Since $m < k$ already holds under Case A above, we consider the following subcases under Case B above:

Case B.1:  $\sqrt{k} < \sigma(\sqrt{k}) < m < m + 1 < k < \sigma(k)$

Under this subcase,
$$m < k \iff \frac{\sigma(k)}{m} > 1 \iff D(k) > 1 \iff k \text{ is not an odd almost perfect number }.$$
Since we ultimately want to prove $m < k$, this case is OK.

Case B.2:  $\sqrt{k} < \sigma(\sqrt{k}) < k < \sigma(k) < m < m + 1$

Note that, under this subcase,
$$0 \leq 2k - \sigma(k) = D(k) = \frac{\sigma(k)}{m} < 1,$$
forcing $k$ to be perfect.  This contradicts the fact that $k$ is a square.
\end{remark}

\section{Conclusion}
We therefore conclude that:

\begin{theorem}\label{thm:2}
Let $n = km$ be a spoof odd perfect number.  Then $m < k$, which is true if and only if $k$ is not an odd almost perfect number.
\end{theorem}

\begin{remark}\label{rem:9}
Note that Theorem \ref{thm:2} is an analogue of a similar result proved by Dris [\cite{4}, \cite{5}] for the case of odd perfect numbers, i.e. the Euler factor $Q^K$ of an odd perfect number $Q^K N^2$ is less than the non-Euler part $N^2$.
\end{remark}

\begin{remark}\label{rem:10}
Virtually everything written in this article (for spoofs) also apply to the case of the usual odd perfect numbers, except for the fact that in the case of spoofs, the quasi-Euler prime $m$ is tacitly assumed to have exponent $1$.  (In the language of the usual odd perfect numbers, the Descartes-Frenicle-Sorli conjecture is a given for spoofs.)
\end{remark}

\section{Acknowledgments}
The author thanks Professor Pace Nielsen for hinting over MathOverflow that most non-computational results on odd perfect numbers carry over to the case of spoofs, and vice-versa.  (Non-computational results means findins on odd perfect numbers that do not depend on numerical computations involving prime factorizations.)  The author also thanks the anonymous referee(s) for valuable feedback that helped in improving the style and presentation of the manuscript.

\end{document}